\documentclass[12pt,a4paper]{report}
\usepackage[cp1251]{inputenc}
\usepackage[russian]{babel}
\usepackage{amsmath,amsthm,layout}
\usepackage{amsfonts}
\usepackage{amssymb}
\usepackage{graphicx}
\usepackage{epstopdf}
\usepackage{wrapfig}
\usepackage{verbatim}
\usepackage{float}
\usepackage{indentfirst, cmap}

\newtheorem{teor}{Теорема}

\newcommand{\R}{\mathbb{R}}

\def\R{{\mathbb R}}

\def\to{\rightarrow}

\renewcommand{\le}{\leqslant}
\renewcommand{\ge}{\geqslant}
\renewcommand{\leq}{\leqslant}
\renewcommand{\geq}{\geqslant}

\renewcommand{\P}{\mathbf P}

\begin{document}
\thispagestyle{empty}  
\large{ \hfuzz=3pt

\begin{titlepage}

\vspace{5mm}
\begin{flushright}
\emph{На правах рукописи}
\end{flushright}

\vspace{20mm}

\begin{centering}

{\bf Кокоткин Андрей Александрович}

\vspace{10mm}

{\bf ВЕРОЯТНОСТНЫЙ ПОДХОД К ЗАДАЧАМ О ГРАФАХ РАССТОЯНИЙ И ГРАФАХ ДИАМЕТРОВ}

\vspace{10mm}


01.01.09 --- дискретная математика и математическая кибернетика

\vspace*{\fill}

{\bf Автореферат} \\
диссертации на соискание ученой степени \\
кандидата физико-математических наук

\vspace*{\fill}

{\bf Москва --- 2014}

\end{centering}
\end{titlepage}

\clearpage
\thispagestyle{empty}

\noindent
Работа выполнена на кафедре анализа данных факультета инноваций и высоких технологий Федерального государственного автономного образовательного учреждения высшего профессионального образования ``Московский физико-технический институт (государственный университет)''.

\vspace{2mm}

\noindent
{\bf Научный руководитель: }
доктор физико-математических наук,
профессор Райгородский Андрей Михайлович.
Место работы: Федеральное государственное бюджетное образовательное учреждение высшего профессионального образования~``Московский государственный университет имени М.В. Ломоносова''.

\noindent
{\bf Официальные оппоненты: }
\begin{itemize}
\item
доктор физико-математических наук, профессор, зав. кафедрой теории вероятностей и дискретной математики ФГБОУ ВПО~``Иркутский государственный университет''  Кузьмин Олег Викторович;
\item
кандидат физико-математических наук, доцент кафедры математических и естественно-научных дисциплин ФГБОУ ВПО~``Российский государственный университет туризма и сервиса'' Лавренченко Сергей Александрович.
\end{itemize}

\noindent
{\bf Ведущая организация:}
Хабаровское отделение Федерального
государственного бюджетного учреждения
науки Института прикладной математики
Дальневосточного отделения Российской академии наук.

\vspace{2mm}

\noindent
Защита состоится 18 декабря 2014 года в 13:00 на заседании диссертационного совета Д 002.017.02 при Федеральном государственном бюджетном
учреждении науки «Вычислительный центр имени А.А. Дородницына Российской академии наук» по адресу 119333, г. Москва, ул. Вавилова, д. 40, конференц-зал.

\vspace{2mm}

\noindent
С диссертацией можно ознакомиться в библиотеке ВЦ РАН и на сайте http://www.ccas.ru/.

\vspace{2mm}

\noindent
Автореферат разослан 11 ноября 2014 года.

\vspace*{\fill}

\noindent
\begin{tabular}{lr}
Ученый секретарь диссертационного совета &\\
доктор физико-математических наук, & \\
профессор & В.~В.~Рязанов \\
\end{tabular}

\clearpage
\newpage
\setcounter{page}{1}
 \pagestyle{plain} \begin{center}
{\Large {\bf Общая характеристика работы}}
\end{center}
{\bf
Актуальность работы}\vspace{0.3em}

Настоящая работа стоит на стыке двух дисциплин: вероятностной комбинаторики и дискретной геометрии. Ниже мы расскажем о ключевых задачах каждой из них.

Дискретная геометрия как наука оформилась на рубеже XIX--XXвв. Одним из ее основоположников можно считать Минковского\footnote{H. Minkowski, \emph{Geometrie der Zahlen}, RG Teubner, Leipzig--Berlin, 1953.}, который исследовал расположение целочисленных векторов по отношению к выпуклым телам в пространстве и доказал фундаментальную теорему о выпуклом теле. Столь же значимые результаты в этой теме получили Вороной\footnote{Г. Вороной  \emph{Собрание сочинений в 3-х томах}, 1952.}, Коркин, Золотарев\footnote{A. Korkine, G. Zolotareff, \emph{Sur les formes quadratiques positives}, Math. Ann. 11 (2), 1877, 242--292.} и другие. Сейчас это отдельный раздел теории чисел и геометрии, называемый геометрией чисел\footnote{Дж. Касселс,  \emph{Геометрия чисел}, Москва, Мир, 1965.}.

С геометрией чисел тесно связано еще одно направление дискретной геометрии. К этому направлению прежде всего относятся следующие три задачи. Первая --- классическая задача Ньютона о плотнейшей упаковке шаров в пространстве\footnote{Л. Тот, \emph{Расположения на плоскостии, сфере и в пространстве}, Москва, Физмалит, 1958.}$^,$\footnote{К. Роджерс, \emph{Укладки и покрытия}, Москва,  Мир, 1968.}$^,$\footnote{Дж. Конвей, Н. Слоэн, \emph{Упаковки шаров, решетки и группы}, Москва, Мир, 1990.}. Вторая задача, двойственная к первой, это задача о редчайшем покрытии пространств шарами. Наконец, третья задача, о контактном числе --- наибольшем количестве шаров, касающихся данного шара в пространстве\footnote{K. Schutte, B.L. van der Waerden,  \emph{Das Problem der dreizehn Kugeln}. Math. Ann. 125, 1953, 325--334.}$^,$\footnote{О. Мусин, \emph{Проблема двадцати пяти сфер}, УМН, 58 (352), 2003, 153--154.}.

Еще одно направление исследований в дискретной геометрии было инициировано Клейн, Эрдешем и Секерешем в 1934 году, которые доказали существование такого числа $N(n)$, что среди любых $N$ точек общего положения на плоскости найдется выпуклый $n$-угольник\footnote{P. Erd\H{o}s, G. Szekeres, \emph{A combinatorial problem in geometry}, Compositio Math., 2, 1935, 463--470.}$^,$\footnote{Harborth, Heiko, \emph{Konvexe Funfecke in ebenen Punktmengen}, Elem. Math. Т. 33 (5), 1978, 116--118.}$^,$\footnote{J. Horton, \emph{Sets with no empty convex 7-gons}, Canadian Math. Bulletin Т. 26 (4), 1983, 482--484.}$^,$\footnote{M. Overmars,  \emph{Finding sets of points without empty convex 6-gons}, Discrete and Computational Geometry Т. 29 (1), 2003, 153--158.}$^,$\footnote{В. Кошелев \emph{Задача Эрдеша--Секереша о пустых шестиугольниках на плоскости}, Моделирование и анализ информационных систем, 16 (2), 2009, 22--74.}.

Следующий класс проблем дискретной геометрии также предложен Эрдешем. Первый вопрос, который был им поставлен, --- о наибольшем числе единичных расстояний в множестве из $n$ точек на плоскости и в пространстве\footnote{P. Erd\H{o}s, \emph{On a set of distances of n points}, Amer. Math. Monthly, 53, 1946, 248--250.}. Второй, не менее важный, наоборот, о числе различных расстояний в таком $n$-точечном множестве\footnote{P. Brass, W. Moser, J. Pach, \emph{Research problems in discrete geometry}, Springer, New York, 2005.}$^,$\footnote{A.M. Raigorodskii, \emph{Coloring Distance Graphs and Graphs of Diameters}, Thirty Essays on Geometric Graph Theory, J. Pach ed., Springer, 2013,
429 - 460.}.

Все перечисленные выше задачи исключительно важны для дискретной геометрии, однако необходимо выделить еще две проблемы, которые имеют особенное значение для нашей работы и также носят фундаментальный характер. Первая из них это гипотеза Борсука\footnote{K. Borsuk, \emph{Drei S\"{a}tze \"{u}ber die $n$-dimensionale euklidische Sph\"{a}re}, Fundamenta Math., 20, 1933, 177--190.} о разбиении множества на части меньшего диаметра. Вторая --- задача Нелсона--Эрдеша--Хадвигера\footnote{P.K. Agarwal, J. Pach, \emph{Combinatorial Geometry}, Wiley-Interscience Series in Discrete Mathematics and Optimization, John Wiley and Sons, 1995.} о хроматическом числе метрических пространств, то есть о наименьшем количестве цветов, в которые так можно его покрасить, чтобы никакие две точки одного цвета не находились на расстоянии единица\footnote{A.M. Raigorodskii, {\it Three lectures on the Borsuk partition problem}, London Mathematical Society Lecture Note Series, 347 (2007), 202--248.}.

Теперь очертим круг задач, характерных для вероятностной комбинаторики. Основной метод, используемый здесь, это вероятностный метод. Его основоположником считается Эрдеш\footnote{P. Erd\H{o}s,  \emph{Some Remarks on the Theory of Graphs}, Bulletin of the American Mathematical Society, 53, 1947, 292--294.}. Вероятностный метод позволил доказать или опровергнуть ряд гипотез, державшихся долгое время\footnote{П. Эрдёш, Дж. Спенсер, \emph{Вероятностные методы в комбинаторике}
Издательство Мир, 1976.}$^,$\footnote{Н. Алон, Дж. Спенсер, \emph{Вероятностный метод} Бином. Лаборатория знаний, 2007.}.

Одним из основных объектов и одновременно инструментов исследования становится случайный граф. В 1959 году Эрдеш и Реньи\footnote{P. Erd\H{o}s, A. R\'enyi, {\it On random graphs, I}, Publ. Math. Debrecen, 6, 1959, 290--297.} дают первую модель случайного графа. Сразу же возникает огромное количество задач об исследовании самых разных характеристик случайного графа, о критических вероятностях для всевозможных его свойств, о распределении какой-либо случайной величины на этом графе. Наша работа посвящена двум свойствам: ``быть графом расстояний'' и ``быть графом диаметров''.

В настоящее время теория графов переживает очередной расцвет в связи с применением ее к исследованиям социальных и биологических и информационных сетей. Теперь, когда стало возможным эмпирически вычислить некоторые характеристики графа, описывающего сеть Интернет, и графов других социальных сетей, выяснилось, что модель Эрдеша--Реньи ни при каких параметрах (а он в этой модели всего один --- вероятность ребра) не может служить сколько-нибудь адекватной моделью для этих графов. В 1999 году Барабаши и Альберт\footnote{A.--L. Barab\'asi, R. Albert,  \emph{Emergence of scaling in random networks} Science, 286, 1999, 509--512.} предложили новую модель случайного графа, а Боллобаш\footnote{B. Bollob\'as, \emph{Random Graphs}, Cambridge Univ. Press, Second Edition, 2001.} и Риордан\footnote{B. Bollob\'as, O.Riordan \emph{The degree sequence of a scale-free random graph process} Random Structures and Algorithms, 18(3), 2001, 279--290.} уточнили ее и доказали, что в ней выполняется по крайней мере два ключевых свойства веб-графа: степенной закон распределения вершин и малый диаметр. С тех пора эта тема получила широкое развитие, были предприняты как попытки улучшения этой модели, так и создания кардинально новых\footnote{А.М. Райгородский, \emph{Модели интернета}, Интеллект, Долгопрудный, 2013.}.

\sloppy

\vspace{0.3em}{\bf Цель работы }

Цель работы состоит в исследовании двух различных свойств случайного графа в модели Эрдеша--Реньи. Во-первых, изучается свойство быть изоморфным графу расстояний; во-вторых, --- свойство быть изоморфным графу диаметров.


 \vspace{0.3em}{\bf Научная новизна}\vspace{0.2em}

Все результаты диссертации являются новыми. Перечислим основные из
них:
\begin{enumerate}
    \item  Доказано, что в любом дистанционном графе на плоскости есть большой индуцированный подграф с хроматическим числом не больше четырех.
    \item  Получены верхние и нижние оценки пороговой вероятности для реализации случайного графа графом расстояний в размерностях от одного до восьми.
    \item  Получены верхние и нижние оценки для размера наибольшего такого подграфа случайного графа, что он представляется как граф диаметров в пространстве произвольной фиксированной размерности.

\end{enumerate}

\vspace{0.3em}{\bf Основные методы исследования}\vspace{0.2em}

В работе используются различные методы теории вероятностей и комбинаторики: методы теории случайных графов, приложения теории плотнейших упаковок, вероятностная техника получения оценки для реализации графом диаметров, связанная с применением теоремы о взаимном расположении подмножеств конечного множества (теорема Р\"{e}для).

\vspace{0.3em}{\bf Теоретическая и практическая
ценность}\vspace{0.2em}

Диссертация имеет теоретический характер. Полученные в диссертации результаты представляют интерес для специалистов в области теории графов, комбинаторики, комбинаторной и дискретной геометрии.

\vspace{0.3em}{\bf Апробация работы}\vspace{0.2em}

Результаты диссертации докладывались на следующих научно-исследовательских конференциях и семинарах:
\begin{itemize}
\item[---] на кафедральном семинаре кафедр дискретной математики и анализа данных (ФИВТ, МФТИ, 2012--2014);

\item[---]  на семинаре ``Вероятностные и алгебраические методы в комбинаторике'' под руководством профессора А.М. Райгородского (механико-математический факультет МГУ, 2009--2012);

\item[---] на международной конференции ``Fete of combinatorics and computer science'' (Кестхей, Венгрия, 11--15 августа 2008);

\item[---] на международной конференции ``Четвертая польская конференция по комбинаторике'' (Бедлево, Польша, 17--21 сентября 2012).
\end{itemize}


\vspace{0.3em}{\bf Публикации}\vspace{0.2em}
\nopagebreak\par
Результаты диссертации опубликованы в 5 работах автора (4 из которых
входят в перечень ВАК), список которых приведен в конце
автореферата. 

\vspace{0.3em}{\bf Структура диссертации}\nopagebreak\vspace{0.2em}

Диссертация состоит из введения, трех глав, заключения и списка литературы,
насчитывающего 68 наименований. Общий объем диссертации составляет
69 страниц.

\begin{center}
{\Large {\bf Краткое содержание диссертации}}
\end{center}

Во введении объясняется актуальность работы и излагается история развития комбинаторной и дискретной геометрии, а также вероятностного метода в комбинаторике. Здесь же очерчен круг ключевых задач этих направлений. 

\vspace{0.3em}{\bf Содержание главы 1}\vspace{0.2em}

В первой главе обсуждается история задач о хроматическом числе пространства и о дистанционных графах, даются необходимые определения, приводятся формулировки и доказательства полученных результатов.

\textit{Хроматическое число} $\chi(\mathbb{R}^{d})$
пространства $\mathbb{R}^{d}$ --- это наименьшее количество
цветов, в которые можно так покрасить $\mathbb{R}^{d}$, чтобы
среди точек одного цвета не нашлось пары точек на расстоянии
единица, то есть
\begin{multline*}
$$
\chi(\mathbb{R}^{d})=\min\,\{k\in
\mathbb{N}:~\exists
V_1,\ldots,V_k,~\mathbb{R}^{d}=V_1\sqcup\ldots\sqcup V_k,\\
~\forall i,~ \forall \textbf{x},\textbf{y} \in
V_i,~\rho(\textbf{x},\textbf{y})\neq 1\},$$
\end{multline*}
где $\rho$ --- обычное евклидово расстояние.

Легко показать, что для любого $d$ величина $\chi(\mathbb{R}^{d})$
конечна. Проблема отыскания хроматического числа пространства была
впервые поставлена на рубеже 40х--50х годов ХХ века.
Несмотря на значительный интерес, вызванный этой проблемой, она до
сих пор, по существу, остается нерешенной. Конечно,
$\chi(\mathbb{R}^{1})=2$, однако уже для плоскости лучшее, что мы
знаем, это оценка $$4\leq\chi(\mathbb{R}^{2})\leq7.$$ Для
трехмерного пространства мы имеем\footnote{D. Coulson, {\it A 15-colouring of 3-space omitting
distance one}, Discrete mathematics, 256, 2002, 83--90.}$^{,}$\footnote{O. Nechushtan, {\it Note on the space chromatic number},
Discrete Mathematics, 256, 2002, 499--507.}
$$6\leq\chi(\mathbb{R}^{3})\leq15,$$ наконец для
растущей размерности\footnote{А.М. Райгородский, {\it Проблема Борсука и
хроматические числа метрических пространств}, Успехи матем. наук,
56(1), 2001, 107--146.}$^,$\footnote{D. Larman, C. Rogers, {\it The realization
of distances within sets in Euclidean space}, Mathematika, 19, 1972, 1--24.}
$$(1.239\ldots+o(1))^d\leq\chi(\mathbb{R}^{d})\leq(3+o(1))^d.$$

Поставленная задача может быть переформулирована в терминах теории графов. Прежде всего \textit{дистанционным
графом} (или \textit{графом расстояний}) назовем конечный граф $G=(V,E)$, вершины которого суть точки евклидова
пространства, а ребра соединяют только пары точек, отстоящих друг от друга на расстояние единица. Иными словами,
$$V\subset\mathbb{R}^{d},~|V|<\infty,~E\subseteq\{(\textbf{x},\textbf{y})\in V\times V:~
\rho(\textbf{x},\textbf{y})=1\}.$$

Напомним, что \textit{хроматическое число} $\chi(G)$ графа $G=(V,E)$ --- это наименьшее количество цветов, в
которые можно так покрасить его вершины, чтобы вершины одного цвета не соединялись ребром, то есть
\begin{multline*}
$$\chi(G)=\min\,\{k\in \mathbb{N}:~\exists V_1,\ldots,V_k,~V=V_1\sqcup\ldots\sqcup V_k,\\~\forall i,~\forall
\textbf{x},\textbf{y}\in V_i,~(\textbf{x},\textbf{y})\notin E\}.$$
\end{multline*}
П. Эрдеш и Н. де Бр\"ейн\footnote{N.G. de Bruijn and P. Erd\H{o}s, {\it A colour
problem for infinite graphs and a problem in the theory of relations},
Proc. Koninkl. Nederl. Acad. Wet., Ser. A, 54 (5), 1951, 371--373.} фактически доказали, что $\chi(\mathbb{R}^{d})=\max\chi(G)$, где максимум
берется по всем графам расстояний в $\mathbb{R}^{d}$. Таким образом, изучение хроматических чисел графов
расстояний играет исключительную роль при исследовании проблемы отыскания хроматического числа пространства.

В этой главе мы сначала рассматриваем случай евклидовой плоскости. Тот факт, что
$\chi(\mathbb{R}^{2})\leq7$, означает, конечно, что для любого графа расстояний $G=(V,E)$ на плоскости
$\chi(G)\leq7$. Как следствие, $\alpha(G)\geq\frac{|V|}{7}$, коль скоро через $\alpha(G)$ мы обозначили
\textit{число независимости} графа $G$, то есть наибольшее количество его вершин, никакие две из которых не
соединены ребром: $$\alpha(G)=\max\,\{|V'|:~ V'\subset V,~\forall \textbf{x},\textbf{y}\in
V',~(\textbf{x},\textbf{y})\notin E \}.$$

Таким образом, в любом ``двумерном'' \ графе расстояний на $n$ вершинах найдется индуцированный подграф,
имеющий не менее $\frac{n}{7}$ вершин и хроматическое число $1$. Это утверждение допускает ряд нетривиальных
обобщений и уточнений. В работе [2] нам удалось доказать следующие результаты для графов расстояний на плоскости. Эти доказательства приводятся в параграфе 1.3.

\vskip+0.2cm
\noindent {\bf Теоремы 1--4.} {\it
Пусть $ k \in \{1, 2, 3, 4\} $. В любом графе расстояний $G=(V,E)$ на $n$ вершинах найдется такой
индуцированный подграф $G'=(V',E') $, что
$ |V'|\geq\left[\frac{kn}{\varkappa}\right] $ и $ \chi(G') \le k,$
где $\varkappa=4.36\ldots$}
\vskip+0.2cm

Для доказательства этих теорем использовался новый результат, являющийся усилением известной теоремы Лармана--Роджерса\footnote{D. Larman, C. Rogers, {\it The realization
of distances within sets in Euclidean space}, Mathematika, 19, 1972, 1--24.}.

\vskip+0.2cm
\noindent {\bf Теорема 8.} {\it
Пусть в $\mathbb{R}^{d}$ задан некоторый граф расстояний
$G=(V,E)$, $\lvert V \rvert =n$. Предположим, существуют такие
числа $k\in\mathbb{N}$ и \, $\nu_0\in(0,1)$, что для любого
индуцированного подграфа $G'=(V',E')$ c $\lvert V' \rvert
\geq[\nu_0 n]$ выполнено $\chi (G')>k$. Тогда для всякого набора
измеримых множеств $S_1, S_2,\ldots, S_k$ в $\mathbb{R}^{d}$ с
условием, что множество $S=S_1\cup S_2\cup\ldots\cup S_k$ имеет
верхнюю плотность $\nu\geq\nu_0$, верно, что каково бы ни было
$a>0$, найдется $S_i$, на котором реализуется расстояние $a$.
}
\vskip+0.2cm

Наиболее интересным является случай $ k = 4 $. Фактически он означает, что в каждом графе расстояний на
плоскости есть индуцированный подграф, который почти целиком совпадает с исходным графом (содержит не менее
91.7\% его вершин) и допускает раскраску в 4 цвета. Если бы в этом утверждении величину 91.7 удалось заменить на
100, то, ввиду теоремы Эрдеша--де Брёйна, это бы означало, что $\chi(\mathbb{R}^{2})=4$.

Зачастую задачи теории графов допускают нетривиальную интерпретацию в терминах \textit{случайного графа}.
Напомним, что одной из наиболее популярных моделей случайного графа является модель, предложенная Эрдешем и А. Реньи\footnote{B. Bollob\'as, {\it Random Graphs}, Cambridge
Univ. Press, Second Edition, 2001.} на рубеже 50х--60х годов XX века. Речь идет о вероятностном пространстве $G(n,p)=(\Omega_n, \mathcal{B}_n,
\P_n)$. Здесь $$\Omega_n=\{G=(V,E):~|V|=n\} \, - $$ множество всевозможных графов на $n$ вершинах (без петель и
кратных ребер), сигма-алгебра $\mathcal{B}_n$ представляет собой множество всех подмножеств $\Omega_n$, а
$$\P_n(G)=p^{|E|}(1-p)^{C_n^2-|E|}.$$ Иначе говоря, можно считать, что ребра случайного графа появляются независимо
друг от друга с вероятностью $p$. Заметим, что в модели Эрдеша--Реньи величина $p$ может зависеть от $ n $.

Нас интересует, с какой вероятностью случайный граф в модели $G(n,p)$ допускает реализацию
на плоскости в качестве графа расстояний. Как это часто бывает в науке о случайных графах, при одних значениях
$p$ эта вероятность будет стремиться к нулю, а при других --- к единице. Определим некоторую критическую
величину $p$, отвечающую за вышеупомянутый ``фазовый переход'', следующим образом:
$$
p^*(n)=
$$
$$
=\sup\left\{p\in[0,1]:\P_n (G \text{ реализуется в } \R^2 \text{ как граф
расстояний})>\frac{1}{2}\right\}.
$$
Такая величина называется {\it пороговой вероятностью}, и она определена корректно ввиду классических результатов о существовании
пороговых вероятностей для монотонных свойств случайных графов.

В той же статье [2] мы доказали следующие результаты.

\addtocounter{teor}{4}
\begin{teor}\label{t2}
При $p=\frac{c}{n}$, где $c<1$, выполнено
$$ \P_n (G \text{ реализуется на плоскости как граф расстояний})\rightarrow 1,\ n\rightarrow\infty.$$
\end{teor}

\begin{teor}\label{t3}
При $p=\frac{c}{n}$, где $c>t_0=14.797\ldots$, выполнено $$ \P_n (G \text{ реализуется на плоскости как граф
расстояний})\rightarrow 0,\ n\rightarrow\infty.$$
\end{teor}

Теоремы \ref{t2} и \ref{t3}, доказанные в параграфе 1.4, означают, что $\frac{1-\varepsilon}{n}\leq p^*(n) \leq\frac{t_0+\varepsilon}{n}$ со сколь угодно
малым положительным $\varepsilon$ и $ n \ge n_0(\varepsilon) $. Тем самым, мы знаем порядок роста пороговой вероятности.

В параграфе 1.5 мы вводим более общую величину
$$
p_d^*(n)=
$$
$$
=\sup\left\{p\in[0,1]:\P_n (G \text{ реализуется в} ~ {\mathbb R}^d ~ \text{как граф
расстояний})>\frac{1}{2}\right\}.
$$
Нам удалось показать, что для размерностей 3--8 также верна
\addtocounter{teor}{4}
\begin{teor}\label{t11}
Положим
$$
c_3 = 55.272\ldots, ~~ c_4 = 164.528\ldots, ~~ c_5 = 504.285\ldots,
$$
$$
c_6 = 1365.170\ldots, ~~ c_7 = 3624.758\ldots, ~~ c_8 = 8675.785\ldots
$$
Тогда при каждом $ d $ и $ p = \frac{c}{n} $, где $ c > c_d $, выполнено
$$
\P_n (G \text{ реализуется в} ~ {\mathbb R}^d ~ \text{ как граф
расстояний})\rightarrow 0,\ n\rightarrow\infty.
$$
\end{teor}
Иными словами,
$$
\frac{1-\varepsilon}{n} \le p_d^*(n) \le \frac{c_d+\varepsilon}{n},
$$
причем величина $ c_d $ растет как экспонента. В одномерном случае порядок роста другой. В пункте 1.5.2 мы доказали, что имеет место следующая

\begin{teor}\label{t4}
Для любого $ \varepsilon > 0 $ существует такое $ n_0 $, что при всех $ n \ge n_0 $ выполнено
$$
\frac{\sqrt[3]{3}-\varepsilon}{n^{4/3}} \le p_1^*(n) \le \frac{\sqrt[3]{12}+\varepsilon}{n^{4/3}}.
$$
\end{teor}

На самом деле в одномерном случае мы знаем точное значение пороговой вероятности.
А именно, в том же пункте 1.5.2 мы доказали, что
$$
\frac{\sqrt[3]{6 \ln 2} - \varepsilon}{n^{4/3}} \le p_1^*(n) \le \frac{\sqrt[3]{6 \ln 2} + \varepsilon}{n^{4/3}}.
$$

\vspace{0.3em}{\bf Содержание главы 2}\vspace{0.2em}

Вторая глава (как и третья) посвящена исследованию некоторых вероятностных характеристик, связанных с классической проблемой Борсука. Напомним, что эта проблема
состоит в отыскании {\it числа Борсука} --- величины $ f(d) $, равной минимальному количеству частей меньшего диаметра, на которые можно разбить
произвольное множество диаметра 1 в пространстве $ {\mathbb R}^d $:
\begin{multline*}
$$
f(d) = \min \{f: ~ \forall ~ \Omega \subset {\mathbb R}^d, ~ {\rm diam}\,\Omega = 1, ~ \exists ~ \Omega_1, \dots, \Omega_f: ~
\Omega = \Omega_1 \cup \ldots \cup \Omega_f, ~\\ \forall ~ i ~~ {\rm diam}\,\Omega_i < 1\}.
$$
\end{multline*}
Гипотеза Борсука\footnote{K. Borsuk, {\it Drei S\"atze \"uber die
$n$-dimensionale euklidische Sph\"are}, Fundamenta Math., 20, 1933}, предложенная своим автором в 1933 году, состояла в том, что $ f(d) = d+1 $. И ровно шестьдесят лет эта гипотеза
оставалась ни доказанной, ни опровергнутой. Лишь в 1993 году Кан и Калаи\footnote{J. Kahn, G. Kalai, \emph{A counterexample to Borsuk’s conjecture}. Bull. Amer. Math. Soc. (N.S.), 29(1), 1993, 60--62} нашли контрпримеры к гипотезе в размерностях $ d \ge 2015 $.
Сейчас известно, что гипотеза Борсука верна при $ d \le 3 $ и ложна при $ d \ge 64 $\footnote{T. Jenrich, \emph{A 64-dimensional two-distance counterexample to Borsuk's conjecture}, 2013, http://arxiv.org/abs/1308.0206v5.}.

С проблемой Борсука тесно связано понятие {\it графа диаметров}. Назовем графом диаметров в $ {\mathbb R}^d $ любой граф $ G = (V,E) $, вершины которого --- точки
$ {\mathbb R}^d $, а ребра --- пары вершин, отстоящих друг от друга на максимальное в множестве вершин расстояние:
$$
V \subset {\mathbb R}^d, ~~ |V| < \infty, ~~
E = \left\{\{{\bf x},{\bf y}\}: ~  |{\bf x}-{\bf y}| = {\rm diam}~V: = \max_{{\bf x}, {\bf y} \in V} |{\bf x}-{\bf y}|\right\}.
$$
Понятно, что минимальное число частей меньшего диаметра, на которые разбивается $ V $ ({\it число Борсука множества} $ V $),
--- это в точности хроматическое число $ \chi(G) $ графа $ G $. Однако было бы некорректно сказать,
что $ f(d) $ --- это максимум таких хроматических чисел. Дело в том, что равенство хроматического числа числу Борсука справедливо лишь в случае конечных
множеств; для бесконечных множеств равенства, вообще говоря, нет: например, если взять в качестве $ V $ сферу в $ {\mathbb R}^d $, то очевидно, что
хроматическое число ее графа диаметров (являющегося паросочетанием) равно двум, тогда как ее число Борсука есть $ d+1 $ ввиду классической
теоремы Борсука--Улама--Люстерника--Шнирельмана\footnote{J. Matou\v{s}ek,
{\it Using the Borsuk--Ulam theorem}, Universitext,
Springer, Berlin, 2003.}.

На плоскости и в пространстве гипотеза Борсука доказана. Более того, существует достаточно много примеров графов диаметров в
$ {\mathbb R}^2 $ и $ {\mathbb R}^3 $ с хроматическими числами 3 и 4 соответственно. Интересно оценить, стало быть, насколько эти примеры часты или редки.

Положим
\begin{multline*}
$$
u_d(n,p)
= \max\left\{k: ~ {\P}_{n,p}\left(\exists ~ H = (W,F) \subset G: \right. \phantom{\frac{1}{2}} |W| = k, \right.
\\
\left. \left.  ~ H = G|_W, ~ H - \text {граф диаметров в } ~ {\mathbb R}^d, ~ \chi(H) = d+1\right) > \frac{1}{2}\right\},
$$
\end{multline*}
\begin{multline*}
$$
u_d'(n,p)
= \max\left\{k: ~ {\P}_{n,p}\left(\exists ~ H = (W,F) \subset G: \right. \phantom{\frac{1}{2}} |W| = k, \right.
\\
\left. \left. ~ H = G|_W, ~ H - \text {связный граф диаметров в } ~ {\mathbb R}^d, ~ \chi(H) = d+1\right) > \frac{1}{2}\right\}.
$$
\end{multline*}

Иными словами, мы ищем максимальное количество вершин в индуцированном (связном) подграфе случайного графа, который одновременно реализуется графом
диаметров в пространстве и имеет при $ d \le 3 $ наибольшее возможное в этом случае хроматическое число (при $ d \ge 64 $ такая постановка
становится несколько произвольной, но по-прежнему осмысленной). Если для любого $ k $
\begin{multline*}
$$
{\P}_{n,p}\left(\exists ~ H = (W,F) \subset G: ~
|W| = k, ~ H = G|_W, ~ H - \text {граф диаметров в } ~ \right.\\ {\mathbb R}^d,  \left.
~ \chi(H) = d+1\right) \le \frac{1}{2},
$$
\end{multline*}
то полагаем $ u_d(n,p) = 0 $ и аналогично поступаем с величиной $ u_d'(n,p) $.

В этой главе сформулированы и доказаны теоремы, полученные нами в [1, 3], верные для размерностей $d=2$, $d=3$. В том случае, когда теоремы верны для обеих величин $u_d$, и $ u_d'$, мы пишем $u_d^*$.
\begin{teor}
Пусть $ p = o\left(\frac{1}{n}\right) $. Тогда при всех достаточно больших $ n \in {\mathbb N} $
выполнено $ u_2^*(n,p) = 0 $.
\end{teor}

\vskip+0.2cm

\begin{teor}
Пусть $ q: = 1-p = o\left(\frac{1}{n^{1.5}}\right) $. Тогда при всех достаточно больших $ n \in {\mathbb N} $
выполнено $ u_2^*(n,p) = 3 $.
\end{teor}

\vskip+0.2cm

\begin{teor}
Пусть $ q = o\left(\frac{1}{n}\right) $, но при этом $ qn^{1.5} \to \infty $.
Тогда при всех достаточно больших $ n \in {\mathbb N} $
выполнено $ u_2^*(n,p) = 4 $.
\end{teor}

\vskip+0.2cm

\begin{teor}
Положим $ \tau(n) = pn $ и $ \sigma(n) = q \ln n $. Пусть $ \tau(n) $ и $ \sigma(n) $ стремятся
к бесконечности с ростом $ n $. Тогда для любого $ \varepsilon > 0 $ существует такое $ n_0 $, что при $ n \ge n_0 $
выполнено
$$
u_2^*(n,p) \le (2+\varepsilon) \log_{\frac{1}{1-p}} (np).
$$
\end{teor}

\vskip+0.2cm

\begin{teor}
Положим $ \tau(n) = \frac{p\sqrt[4]{n}}{\ln n} $ и $ \sigma(n) = q \ln n $. Пусть $ \tau(n) $ и $ \sigma(n) $ стремятся
к бесконечности с ростом $ n $. Тогда для любого $ \varepsilon > 0 $ существует такое $ n_0 $, что при $ n \ge n_0 $
выполнено
$$
u_2^*(n,p) \ge \left(2-\varepsilon+\frac{4 \ln p}{\ln (np)}\right) \log_{\frac{1}{1-p}} (np).
$$
\end{teor}

\vskip+0.2cm

\begin{teor}
Зафиксируем некоторое число $ \alpha \in \left(0, \frac{1}{2}\right) $ и положим $ \tau(n) = pn^\alpha $.
Пусть с некоторым $ C > 0 $ начиная с некоторого $ n $ выполнено $ 1 < \tau(n) < C $. Тогда для
любого $ \varepsilon > 0 $ существует такое $ n_0 $, что при $ n \ge n_0 $ имеет место неравенство
$$
u_2^*(n,p) \ge \left(2-2\alpha-\varepsilon\right)\frac{\ln n}{p}.
$$
\end{teor}

\vskip+0.2cm

\begin{teor}
Пусть найдется такое $c<1$, что
$p<\frac cn$. Тогда при всех достаточно больших $n\in {\mathbb N}$
выполнено $ u_3^*(n,p) = 0 $.
\end{teor}

\vskip+0.2cm

\begin{teor}
Пусть $ q = o\left(\frac{1}{n^{1.5}}\right) $. Тогда при всех достаточно
больших $ n \in {\mathbb N} $ выполнено $ u_3^*(n,p) = 4 $.
\end{teor}

\vskip+0.2cm

\begin{teor}
Положим $ \tau(n) = pn $ и $
\sigma(n) = q \ln n $. Пусть $ \tau(n) $ и $ \sigma(n) $ стремятся
к бесконечности с ростом $ n $. Тогда для любого $ \varepsilon > 0
$ существует такое $ n_0 $, что при $ n \ge n_0 $ выполнено
$$
u_3^*(n,p) \le (2+\varepsilon) \log_{\frac{1}{1-p}} (np).
$$
\end{teor}

\vskip+0.2cm

\begin{teor}
Пусть для всякого $ \alpha > 0 $
выполнено $ pn^{\alpha} \to \infty $ и $q\ln n \to \infty$ при
$n\to \infty$. Тогда для любого $ \varepsilon > 0 $ существует
такое $ n_0 $, что при $ n \ge n_0 $ выполнено
$$
u_3^*(n,p) \ge \left(2-\varepsilon\right) \log_{\frac{1}{1-p}} (np).
$$
\end{teor}

\vskip+0.2cm

\begin{teor}
Зафиксируем некоторое
$\alpha \in \left(0,\frac14\right)$ и положим $ \tau(n) = pn^\alpha $.  Пусть с некоторым $ C > 0 $ начиная с некоторого $ n $ выполнено $
1<\tau(n)<C $. Тогда для любого $ \varepsilon > 0 $ существует такое $ n_0 $, что при $ n \ge n_0 $ имеет место неравенство
$$
u_3^*(n,p) \ge \left(2-4\alpha-\varepsilon\right)\frac{\ln n}{p}.
$$
\end{teor}

\vskip+0.2cm

Видно, что при разумных условиях на асимптотику вероятности ребра теоремы 16 и 17 дают асимптотику величины $ u_2^*(n,p) $ и эта асимптотика имеет вид
$ 2\log_{\frac{1}{1-p}} (np) $. Более того, если $ p $ стремится к нулю, то в условиях теоремы 17 выполнено $ 2\log_{\frac{1}{1-p}} (np) \sim
2\frac{\ln n}{p} $. Поэтому теорема 18 просто дает аналог оценки из теоремы 17, который лишь в константу раз хуже, но работает на большем диапазоне
значений $ p $. Полностью аналогичная картина имеет место в теоремах 21--23. При этом утверждения теорем 16 и 21 совсем идентичны, в теореме 22 практически
та же оценка, что и в теореме 17, но более узкий диапазон допустимых вероятностей ребра, а в теореме 23 и диапазон уже, чем в теореме 18, и оценка послабее.

В параграфе 2.3 мы доказываем теоремы 13--18, в которых идет речь о величине $u_2^*(n,p)$. Величине $u_3^*(n,p)$ мы посвящаем параграф 2.4, здесь доказаны теоремы 19--23.

\vspace{0.3em}{\bf Содержание главы 3}\vspace{0.2em}

В этой главе сформулированы и доказаны результаты о величинах $u_d(n,p)$ и $u'_d(n,p)$, введенных в главе 2, верные для произвольной фиксированной размерности. Эти результаты мы получили в работе [5].

\vskip+0.2cm

\begin{teor}
Пусть для всякого $ \alpha > 0 $
выполнено $ pn^{\alpha} \to \infty $ и $q\ln n \to \infty$ при
$n\to \infty$. Тогда для любого $ d $ и для любого $ \varepsilon > 0
$ существует такое $ n_0 $, что при $ n \ge n_0 $ выполнено
$$
u_d(n,p) \ge (2-\varepsilon) \log_{\frac{1}{1-p}} (np).
$$
\end{teor}

\vskip+0.2cm

Таким образом, нижняя оценка справедлива в любой размерности, но в чуть разных условиях (для размерности 2 ограничений меньше).
Повторим, однако, что если для размерностей 2 и 3 эта оценка выполнялась
также для величины со штрихом, то при $ d \ge 4 $ метод немного другой и на величину со штрихом он не распространяется.

\vskip+0.2cm

\begin{teor}
Зафиксируем некоторое
$\alpha \in \left(0,\frac{1}{d}\right)$ и положим $ \tau(n) = pn^\alpha $.  Пусть с некоторым $ C > 0 $ начиная с некоторого $ n $ выполнено $
1<\tau(n)<C $. Тогда для любого $ d $ и любого $ \varepsilon > 0 $ существует такое $ n_0 $, что при $ n \ge n_0 $ имеет место неравенство
$$
u_d(n,p) \ge \left(2-2\alpha-\varepsilon\right)\frac{\ln n}{p}.
$$
\end{teor}

\vskip+0.2cm

\begin{teor}
Пусть $ p $ --- это либо константа, либо произвольная функция, которая стремится к нулю при $ n \to \infty $, но при этом
ограничена снизу величиной $ \frac{c}{n} $, где $ c > 1 $. Тогда для любого $ d $ и для любого $ \varepsilon > 0
$ существует такое $ n_0 $, что при $ n \ge n_0 $ выполнено
$$
u_d^*(n,p) \le (2+\varepsilon) (d+1) \log_{\frac{1}{1-p}} (np).
$$
\end{teor}

\vskip+0.2cm

Теорема 26 дает абсолютно универсальную верхнюю оценку, работающую даже в более широком диапазоне, нежели теорема 16. В этом ее большой плюс. Однако при $ d = 2, 3 $
значение полученной в ней оценки несколько хуже ранее известных. Следующая теорема устраняет эту проблему ценой сокращения диапазона допустимых вероятностей ребра.

\vskip+0.2cm

\begin{teor}
Пусть $ p $ --- это константа. Тогда для любого $ d $ и для любого $ \varepsilon > 0
$ существует такое $ n_0 $, что при $ n \ge n_0 $ выполнено
$$
u_d^*(n,p) \le (2+\varepsilon) \left[\frac{d}{2}\right] \log_{\frac{1}{1-p}} (np).
$$
\end{teor}

\vskip+0.2cm

Теорема 27 отлично согласуется с теоремами 16 и 21.
Тем не менее, в ней $ p $ --- константа. Слегка оторваться от этого ограничения можно. Но сделать это довольно тяжело, т.к. доказательство теоремы 27 опирается на ряд тонких утверждений, многие из которых с трудом обобщаются на случаи непостоянных вероятностей. Отметим также, что при $ d \ge 4 $ нижние и верхние оценки начинают отличаться друг от друга, так что асимптотику найти уже не удается.

Наконец, обобщением и даже усилением теорем 13 и 19 служит следующая


\begin{teor}
Пусть $ d \ge 3 $ и найдется такое $c < 2(d-1)\ln (d-1)$, что
$p<\frac cn$. Тогда при всех достаточно больших $n\in {\mathbb N}$
выполнено $ u_d^*(n,p) = 0 $.
\end{teor}

Теоремы 26, 28 доказаны в параграфе 3.2, теоремы 24, 25 --- в параграфе 3.3, наконец, теореме 27 посвящен четвертый параграф. В пятом параграфе рассмотрены сложности обобщения последней теоремы на случай убывающей вероятности ребра $p$.

В заключении даны возможные направления дальнейших исследований.

\nopagebreak

\vspace{0.3em}{\bf Благодарности}\vspace{0.2em}

 Автор признателен профессору Андрею Михайловичу
Райгородскому за постановку задач и неоценимую помощь в
работе.
\begin{center}

\newpage
{\bf Список публикаций по теме диссертации}\\
\end{center}

\begin{enumerate}
\item[{[1]}] А.А. Кокоткин, А.М. Райгородский, {\it О реализации случайных графов графами диаметров}, Труды МФТИ, 4 (2012), N1, 19 -- 28.


\item[{[2]}] А.А. Кокоткин, А.М. Райгородский, {\it О больших подграфах графа расстояний, имеющих маленькое хроматическое число},
Современная математика. Фундаментальные направления, 51 (2013), 64 -- 73.

\item[{[3]}] А.А. Кокоткин, А.М. Райгородский,  {\it О реализации подграфов случайного графа графами диаметров на плоскости и в пространстве}, Труды МФТИ, 6 (2014), N2, 44 -- 60.

\item[{[4]}] А.А. Кокоткин, {\it О реализации подграфов случайного графа графами диаметров в евклидовых пространствах}, Доклады РАН, 456 (2014), N6, 1 -- 3.

\item[{[5]}] А.А. Кокоткин, {\it О больших подграфах графа расстояний, имеющих маленькое хроматическое число}, Математические заметки, 96 (2014), N 2, 318 -- 320.

\end{enumerate}

\appendix

\newpage
\bibliographystyle{srt}
\bibliography{5.Literature/Literature}
\end{document}